\documentclass[a4paper,11pt]{article}
\usepackage{array}
\usepackage{anysize}
\usepackage{url}
\usepackage{paralist}
\usepackage{amsthm}
\usepackage{amsmath}
\usepackage{amssymb}
\usepackage{amscd} 
\usepackage{bbm}
\usepackage{graphicx}
\usepackage{psfrag}
\usepackage{pstricks}
\usepackage{subfigure}

\newcommand\gon{m}
\newcommand\T{\mathsf{T}}
\newcommand\ZO{\{0,1\}}
\newcommand\id[1]{\mathsf{Id}_{#1}}
\newcommand\R{\mathbb{R}}
\newcommand\Hr{\mathcal{H}}
\newcommand\Vr{\mathcal{V}}
\newcommand\aff{\mathsf{aff}\,}
\renewcommand\dim{\mathsf{dim}\,}
\renewcommand\min{\mathsf{min}}
\newcommand\conv{\mathsf{conv}\,}
\newcommand\V{\mathsf{vert}\,}
\newcommand\eps{\varepsilon}
\newcommand\PM{\{+,-\}}
\newcommand\PMZ{\{+,-,0\}}
\newcommand\eq{\mathsf{eq}\,}
\newcommand\s{\sigma}
\newcommand\ovA{\overline{A}}
\newcommand\oovA{\overline{\ovA}}
\newcommand\ovG{\overline{G}}
\newcommand\hog{^\mathsf{hog}}
\newcommand\1{\mathbbmss{1}}
\newcommand\0{\mathbb{O}}
\newcommand\la{\lambda}
\newcommand\ncp{\mathsf{NCP}}
\newcommand\pdpp[2]{\mathsf{PDPP}_{#1}(#2;\gon)}
\newcommand\dpp[2]{P_{#1}(#2;\gon)}
\newcommand\cyc{\mathsf{Cyc}}
\newcommand\w{\omega}
\newcommand\Qw{\tilde{Q}}
\newcommand\ovw{\overline{\w}}
\newcommand\Lex[1]{\mathsf{Lex}_{#1}\,}
\renewcommand\L{\mathcal{L}}
\renewcommand\a{\alpha}
\renewcommand\b{\beta}
\newcommand\ngonComb{\mathcal{P}_{\gon}}

\newcommand\FL{\mathcal{F}}
\newcommand\xp{\overline{x}}
\newcommand\xpp{\overline{\xp}}
\newcommand\tP{\tilde{P}}
\newcommand\tF{\tilde{F}}

\newtheorem{thm}{Theorem}[section]
\newtheorem{cor}[thm]{Corollary}
\newtheorem{lem}[thm]{Lemma}
\newtheorem{prop}[thm]{Proposition}
\theoremstyle{definition}
\newtheorem{dfn}[thm]{Definition}

\begin{document}

\title{Construction and Analysis\\of Projected Deformed Products}

\author{
Raman Sanyal \qquad G\"{u}nter M. Ziegler\\\ \\
\small Institute of Mathematics, MA 6-2\\
\small TU Berlin\\\small D-10623 Berlin, Germany
\\\small\tt\{sanyal,ziegler\}@math.tu-berlin.de
}

\date{October 10, 2007}

\maketitle

\begin{abstract}\noindent
    We introduce a deformed product construction for simple polytopes in terms
    of lower-triangular block matrix representations.  We further show how
    Gale duality can be employed for the construction and for the analysis of
    deformed products such that specified faces (e.g.\ all the $k$-faces) are
    ``strictly preserved'' under projection.

    Thus, starting from an arbitrary neighborly simplicial $(d-2)$-polytope
    $Q$ on $n-1$ vertices we construct a deformed $n$-cube, whose projection
    to the last $d$ coordinates yields a \emph{neighborly cubical
    $d$-polytope}.  As an extension of the cubical case, we construct matrix
    representations of deformed products of (even) polygons (DPPs), which have
    a projection to $d$-space that retains the complete $(\lfloor \tfrac{d}{2}
    \rfloor - 1)$-skeleton.

    In both cases the combinatorial structure of the images under projection
    is completely determined by the neighborly polytope $Q$: Our analysis
    provides explicit combinatorial descriptions. This yields a multitude of
    combinatorially different neighborly cubical polytopes and DPPs.

    As a special case, we obtain simplified descriptions of the neighborly
    cubical polytopes of Joswig \& Ziegler (2000) as well as of the
    \emph{projected deformed products of polygons} that were announced by
    Ziegler (2004), a family of $4$-polytopes whose ``fatness'' gets
    arbitrarily close to~$9$.
\end{abstract}

\section{Introduction}\label{sec:intro}

Some remarkable geometric effects can be achieved for projections of
``suitably-deformed'' high-dimensional simple polytopes. This includes the
Klee-Minty cubes \cite{KlMi}, the Goldfarb cubes~\cite{gold82}, and many other
exponential examples for variants of the simplex algorithm, but also the
``neighborly cubical polytopes'' first constructed by Joswig \& Ziegler
\cite{jos00}.  A geometric framework for ``deformed product'' constructions
was provided by Amenta \& Ziegler \cite{am98}.

Here we introduce a \emph{generalized deformed products} construction. In
terms of this construction, the previous version by Amenta \& Ziegler
concerned deformed products of rank~$1$.  The new construction is presented in
matrix version (that is, as an $\mathcal H$-polytope).  Iterated deformed
products are thus given by lower-triangular block matrices, where the blocks
below the diagonal do not influence the combinatorics of the product (for
suitable right-hand sides).

The deformed products $P$ are constructed in order to provide interesting
images after an affine projection $\pi: P\rightarrow \pi(P)$.  The
deformations we are after are designed so that certain classes of faces of the
deformed product $P$, e.g.\ all the $k$-faces, are ``preserved'' by a
projection to some low-dimensional space, i.e.\ mapped to faces of~$\pi(P)$.
In the combinatorially-convenient situation, the faces in question are
\emph{strictly preserved} by the projection; we give a linear algebra
condition that characterizes the faces that are strictly preserved (Projection
Lemma \ref{lem:projection}). We also identify a situation when all nontrivial
faces of $\pi(P)$ arise as images $\pi(F)$ of faces $F\subset P$ that are
strictly preserved (Corollary~\ref{cor:all_strictly}).

The conditions dictated by the Projection Lemma may be translated via a
non-standard application of Gale duality \cite[Sect.~6.3]{grue03}
\cite[Lect.~6]{zie95} into conditions about the combinatorics of an auxiliary
polytope~$Q$.

As an instance of this set-up, we show how \emph{neighborly cubical
$d$-polytopes} arise from projections of a deformed $n$-cube where all the
$(\lfloor\tfrac d2\rfloor-1)$-faces are preserved by the projection. The
precise form of the matrix representation of the $n$-cube, and the
combinatorics of the resulting polytopes, is dictated via Gale duality by a
neighborly simplicial\,(!) $(d-2)$-polytope with $n-1$ vertices.  As special
cases, we obtain the neighborly cubical polytopes first obtained by Joswig \&
Ziegler \cite{jos00}, and also geometric realizations for neighborly cubical
spheres as described by Joswig \& R\"orig~\cite{jos05}.

Finally, we construct and analyze \emph{projected deformed products of (even)
polygons} (PDPP polytopes), as the images of a deformed product of $r$ even
polygons, projected to~$\R^d$. The projection is designed to strictly preserve
all the $(\lfloor\tfrac{d}{2}\rfloor-1)$-faces (as well as additional $\tfrac
d2$-faces, if $d$ is even).  This produces in particular the $2$-parameter
family of $4$-dimensional polytopes from \cite{zie04}, for which the
``fatness'' parameter introduced in \cite{zie02} gets as large as
$9-\varepsilon$.  We present a new construction  (drastically simplified and
systematized) and a complete combinatorial description of these polytopes.

This work is based on the Diploma Thesis \cite{dipl-Sanyal}; see also the
research announcements in  \cite{zie04} and \cite{Z102a}.  The ``wedge
product'' polytopes of R\"orig \& Ziegler \cite{Z109} provide another
interesting instance of ``deformed high-dimensional simple polytopes''.  A
further analysis shows that the neighborly cubical polytopes, the PDPP
polytopes as well as the wedge products do exhibit a wealth of interesting
polyhedral surfaces, including the ``surfaces of unusually high genus'' by
McMullen, Schulz \& Wills \cite{mcmullen83:_polyh_e}, and equivelar surfaces
of type $(p,2q)$.  Topological obstructions that prevent a suitable projection
of ``deformed products of \emph{odd} polygons'', or of the wedge product
polytopes, will be presented by R\"orig \& Sanyal \cite{roe07}.

{\bf Acknowledgements.} The first author would like to thank Andreas
Paffenholz, Thilo R\"orig, Jakob Uszkoreit, Arnold Wa{\ss}mer, and Axel Werner
for ``actively listening'' and Vanessa K\"a\"ab for more.  Both authors
gratefully acknowledge support by the German Science Foundation DFG via the
Research Training Group ``Methods for Discrete Structures'' and a Leibniz
grant.

\section{Basics}\label{sec:basics}

In this section we recall basic properties and notation about the main objects
of this paper: \emph{convex polytopes}. Readers new to the country of
polytopia will find useful information in the well-known travel guides
\cite{grue03} and \cite{zie95} while the frequent visitors might wish to skim
the section for possibly non-standard notation.

One of the main messages this article tries to convey is that it pays off to
work with polytopes \emph{in explicit coordinates} (matrix representation).
Classically, there are two fundamental ways of viewing a polytope in
coordinates: the interior or $\Vr$-representation, and the exterior or
$\Hr$-representation.  For $\Vr$-polytopes ``with few vertices'', Perles
\cite[Chap.~6]{grue03} had developed Gale duality as a powerful tool. In this
article, we will apply Gale duality for the analysis of projected simple
$\Hr$-polytopes.  The basics for this will be developed in this section.

\subsection{Polytopes in coordinates}

For the rest of the section, let $P \subset \R^d$ be a full-dimensional
polytope.  In its \emph{interior} or $\Vr$\emph{-presentation}, $P = \conv V$
is given as the convex hull of a finite point set $V  = \{ v_1, \dots, v_m \}
\subset \R^d$ and $V$ is inclusion-minimal with respect to this property. The
elements of $V$ are called \emph{vertices}, with notation $\V P = V$. For a
nonempty subset $I \subseteq [m] = \{1,\dots,m\}$ the set $V_I = \{v_i : i \in
I\}$ forms a \emph{face} of $P$ if there is a linear functional $\ell : \R^d
\rightarrow \R$ such that $\ell$ attains its maximum over~$P$ on $F = \conv
V^\prime$. The dimension $\dim F$ is the dimension of its affine span.  The
empty set is also a face of~$P$, of dimension~$-1$. The collection $\FL P$ of
all faces of $P$, ordered by inclusion, is a graded, atomic and coatomic
lattice with $\dim+1$ as its rank function.  We denote by $\FL\partial P :=
\FL P \setminus \{ P \}$ the face poset of the boundary of $P$.  We say that
two polytopes are of the same \emph{combinatorial type} if their face lattices
are isomorphic as abstract posets. A polytope $P$ is \emph{simplicial} if
small perturbations applied to the vertices do not alter its combinatorial
type.  Equivalently, every $k$-face of $P$ ($k<\dim P$) is the convex hull of
exactly $k+1$ vertices. The \emph{quotient} $P / F$ of $P$ by a face $F$ is a
polytope with face lattice isomorphic to $\FL P_{\ge F} = \{ G \in \FL P : F
\subseteq G \}$.  If $F = \{v \}$ is a vertex, then $P/v$ is called a
\emph{vertex figure} at $v$. 

The polytope $P$ is given in its \emph{exterior} or $\Hr$\emph{-presentation}
if $P$ is the intersection of finitely many halfspaces. That is, if there are
(outer) normals $a_1,\dots,a_n \in \R^d$ and displacements $b_1,\dots,b_n \in
\R$ such that
\[
    P\ =\ \bigcap_{i=1}^n \{ x \in \R^d : a_i^\T x \le b_i \},
\]
where we assume that the collection of normals is \emph{irredundant}, thus
discarding any one of the halfspaces changes the polytope.  The hyperplanes
$H_i = \{ x \in \R^d : a_i^\T x = b_i \}$ are said to be \emph{facet
defining}; the corresponding $(d-1)$-faces $F_i = P \cap H_i$ are called
\emph{facets}.  More compactly, we think of the normals $a_i$ as the rows of a
matrix $A \in \R^{n \times d}$ and, with $b \in \R^n$ accordingly, write 
\[
    P = P(A,b) = \{ x \in \R^d : A\,x \le b \}.
\]
For any subset $F \subseteq P$ let $\eq F = \{ i \in [n] : F \subset H_i \}
\subseteq [n]$ be its \emph{equality set}.  Clearly, $F \subseteq 
\bigcap_{i\in\eq F} F_i$; in case of equality, the set $F$ is a face of $P$.
Denote by $A_I$ the submatrix of $A$ induced by the row indices in $I
\subseteq [n]$. Thus any face $F$ is given by $F = P \cap \{ x : A_I \, x =
b_I \}$, for $I=\eq F$.  The collection of equality sets of faces ordered by
reverse inclusion is isomorphic to $\FL P$.  The polytope $P$ is \emph{simple}
if its combinatorial type is stable under small perturbations applied to the
bounding hyperplanes.  Equivalently, every nonempty face $F$ is contained in
no more than $|\eq F| = d - \dim F$ facets.

\subsection{Gale duality}

Let $P \subset \R^d$ be a $d$-polytope and let the rows of $V \in \R^{m \times
d}$ be the $m$ vertices of $P$. Denote by $V\hog = (V, \1) \in \R^{m \times
(d+1)}$ the homogenization of $V$. The column span of $V$ is a $d+1$
dimensional linear subspace. Choose $G \in \R^{m \times (m-d-1)}$ such that
the columns form a basis for the orthogonal complement. Any such basis,
regarded as an ordered collection of $m$ row vectors, is called a \emph{Gale
transform} of $P$. It is unique up to linear isomorphism and, by the reverse
process, characterizes $V\hog$, again up to linear isomorphism. So it
determines $P$ only up to a projective transformation.  However, the striking
feature of Gale transforms is that its combinatorial properties are, in a
precise sense, dual to those of $P$; this correspondence goes by the name of
\emph{Gale duality}. 

In order to state and work with Gale duality we introduce some concepts and
notations.  As before we write $V_I$ for the subset of the rows of $V$ indexed
by $I \subseteq [m]$. A subset $I \subset [m]$ names a \emph{coface} of $P$ if
the complement $V_{[m] \setminus I}$ is the vertex set of a face of $P$.

\begin{dfn}
    A collection of vectors $G = \{ g_1, g_2, \dots, g_m \} \subset \R^k$ is
    \emph{positively dependent} if there are numbers $\la_1, \la_2, \dots,
    \la_m > 0$ such that $\la_1 g_1 + \cdots + \la_m g_m = 0$.  It is
    \emph{positively spanning} if in addition $G$ is of full rank~$k$.
\end{dfn}

``Begin positively spanning'' is, like ``being spanning'', an \emph{open}
condition, i.e.\ preserved under (sufficiently small) perturbations of the
elements of $G$. This, however, is not true for ``being positively
dependent'': Consider e.g.\ $\{g, -g \}$ for $g \in \R^k$, $g\neq0$, $k>1$.

\begin{thm}[Gale duality]
    Let $P = \conv V$ be a polytope and $G$ a Gale transform of $P$. Let $I
    \subset [m]$ then $I$ names a coface of $P$ if and only if $G_I$ is
    positively dependent.
\end{thm}

In light of Gale duality, the preceding theorem implies that for a general
polytope not every subset of the vertex set of a face necessarily forms a
face. This, however, is true for simplicial polytopes and, in fact,
characterizes them.  A still stronger condition is satisfied if no $d+1$
vertices of a $d$-polytope lie on a hyperplane, that is, if the vertices are
\emph{in general position} with respect to affine hyperplanes.  (The polytope
is then automatically simplicial; however, the vertices of a regular
octahedron are not in general position). This translates into Gale diagrams as
follows.

\begin{prop}\label{prop:simplicial_gale}
    Let $P\subset \R^d$ be a polytope and $G\subset\R^k$ a Gale transform of
    $P$. Then $P$ is simplicial with vertices in general position if and only
    if the rows of $G$ are in general position with respect to linear
    hyperplanes, that is, if any $k$ vectors of $G$ are linearly independent.
\end{prop}

\subsection{Faces strictly preserved by a projection}

Projections are fundamental in polytope theory: Every polytope on $n$ vertices
is the image of an $(n-1)$-simplex under an affine projection. This in
particular says that the analysis of the images of polytopes under projection
is as difficult as the general classification of all combinatorial types of
polytopes.  The problem is that a $k$-face $F\subset P$ can behave in various
ways under projection: It can map to a $k$-face, or to part of a $k$-face, or
to a lower-dimensional face of $\pi(P)$.  Even if it maps to a $k$-face
$\pi(F)\subset\pi(P)$, there may be other $k$-faces of $P$ that map to the
same face $\tF = \pi(F)$.  In that case, the face $\pi^{-1}(\tF)$ has higher
dimension than~$F$.  Thus, as a serious simplifying measure, we restrict our
attention in the following to the most convenient situation, of faces that are
``strictly preserved'' by a projection.

\begin{dfn}[Strictly preserved faces \cite{zie04}]
    \label{dfn:strictly_preserved}
    Let $P$ be a polytope and $Q = \pi(P)$ the image of $P$ under an affine
    projection $\pi: P \rightarrow \R^d$. A nonempty face $F$ of $P$ is
    \emph{(strictly) preserved} by $\pi$ if%
    \begin{compactenum}[(i)]
        \item $\pi(F)$ is a face of $Q$ combinatorially equivalent to $F$, and
            \hfill (preserved face)
        \item the preimage $\pi^{-1}(\pi(F))$ is $F$. 
            \hfill (strictly preserved)
    \end{compactenum}
\end{dfn}

Since in the following we will be concerned exclusively with the analysis of
strictly preserved faces, we will generally drop the modifier ``strictly''
starting now.

The following lemma gives an algebraic way to \emph{read off} the preserved
faces from a polytope in exterior presentation.   Every affine projection
$\pi: \R^n \rightarrow \R^d$ factors as an affine transformation followed a
projection $\pi_d : \R^{n-d} \times \R^d \rightarrow \R^d$ that deletes the
first $n-d$ coordinates, that is $\pi_d(\xp,\xpp) = \xpp$ for all $(\xp,\xpp)
\in \R^{n-d} \times \R^d$.  Therefore, we will focus on the projections
$\pi_d$ ``to the last $d$ coordinates''.  For a polytope $P = P(A,b) \subset
\R^n$ in exterior presentation the projection map $\pi_d$ naturally partitions
the columns of $A$, as $A = (\ovA | \oovA )$.

\begin{lem}[Projection Lemma: Matrix version]\label{lem:projection}
    Let $P = P(A,b) \subset \R^n$ be a polytope, $F$ a nonempty face of $P$,
    and $I = \eq F $ the index set of the inequalities that are tight at~$F$.
    Then $F$ is preserved by the projection $\pi_d:P\rightarrow\R^d$ to the
    last $d$ coordinates if and only if the rows of $\ovA_I$ are positively
    spanning.
\end{lem}

The proof makes use of the following geometric version of the Farkas Lemma.

\begin{lem}[{\cite[Sect.~1.4]{zie95}}] \label{lem:farkas}
    Let $P = P(A,b)$ be a polytope and $F \subseteq P$ a nonempty face. For a
    linear functional $\ell(x)=cx$ we denote by $P^\ell$ the nonempty face of
    $P$ on which $\ell$ attains its maximum.  The linear function $\ell$
    singles out $F$, that is $P^\ell = F$, if and only if $c$ is a strictly
    positive linear combination of the rows of $A_{\eq F}$.
\end{lem}

\begin{proof}[Proof of Lemma \ref{lem:projection}]
    We split the proof into two parts. 

    {\bf Claim 1.} $\tF = \pi_d(F)$ is a face of $\tilde{P} = \pi_d(P)$ with
    $\pi_d^{-1}(\tF) \cap P = F$ iff $\ovA_I$ is positively dependent.

    By Lemma~\ref{lem:farkas} the rows of $\ovA_I$ are positively dependent if
    and only if there is some $c \in \R^d$ such that the linear function
    $\ell(x) := (0,c)\,x = c\, \xpp$ satisfies $P^\ell = F$. Rewriting $\ell =
    h \circ \pi_d$ with $h(\xpp) = c\, \xpp$ we see that such a $c$ exists if
    and only if there is a linear function $h$ on $\tP$ such that $\tP^h =
    \tF$.

    {\bf Claim 2.} Considering $F$ as a (sub-)polytope in its own right, then
    $\tF = \pi_d(F)$ is combinatorially equivalent to $F$ if and only if $A_I$
    has full row rank.

    The polytopes $F$ and $\tF$ are combinatorially equivalent iff they are
    affinely isomorphic. This happens if and only if the linear map $\pi_d$ is
    injective restricted the linear space $L = \{x : A_I\,x = 0\}$, which is
    parallel to $\aff F = \{ x : A_I\,x = b_I\}$ the affine hull of $F$. Now,
    $\pi_d|_L$ is injective iff $\mathsf{ker}\, \pi_d \cap L \cong \{ \xp :
    \ovA_I\, \xp = 0 \}$ is trivial.
\end{proof}

See \cite{san07} for a proof in a different wording.

Lemma~\ref{lem:projection} allows us to guarantee that in certain situations
every single $k$-face is preserved by a projection $\pi:P\rightarrow\pi(P)$.
Then, however, we want to also see that $\pi(P)$ has no other $k$-face than
those induced by the projection. This will be argued via the following lemma.

\begin{lem}\label{lem:all_strictly}
    Let $P = P(A,b) \subset \R^n$ be an $n$-polytope such that for every
    vertex $v \in \V P$ the rows of the matrix $\ovA_{\eq v}$ are in general
    position with respect to linear hyperplanes. Then every proper face of $P$
    is either preserved under $\pi_d$ or falls short of being a face of
    $\pi_d(P)$.
\end{lem}

\begin{proof}
    If $G \subset \R^k$  is a set of at least $k$ vectors in general position
    with respect  to linear hyperplanes then $\dim \mathsf{span}\,G^\prime \ge
    \min\{|G^\prime|,k\}$ for every subset $G^\prime \subseteq G$.  In
    particular, every positively dependent subset is positively spanning.

    Let $F \subset P$ be a proper face. From the proof of Lemma
    \ref{lem:projection} it follows that $\pi_d(F)$ is a face iff $\ovA_{\eq
    F}$ is positively dependent. Let $v \in \V P$ be a vertex with $v \in F$.
    Then $\ovA_{\eq F} \subseteq \ovA_{\eq v}$ and $\ovA_{\eq v} \subset
    \R^{n-d}$ is a set of at least $n$ vectors in general position with
    respect to linear hyperplanes.
\end{proof}

\begin{cor}\label{cor:all_strictly}  
    If \emph{all} $k$-faces of $P$ are preserved by the projection $\pi: P
    \rightarrow \pi(P)$, then all $k$-faces of~$\pi(P)$ arise as images of
    $k$-faces of~$P$.
\end{cor}

\begin{proof}
    For any $k$-face  $G\subseteq \pi(P)$ we know that $\widehat G =
    \pi_d{}^{-1}(G)$ is a face of~$P$, of dimension $\dim \widehat G\ge k$.
    Now if $F\subseteq \widehat G$ is any $k$-face of $\widehat G$, then by
    Lemma~\ref{lem:all_strictly} either $F$ is preserved, and we get
    $\pi_d(F)=G$, or $F$ is not mapped to a face. The latter case cannot arise
    here.
\end{proof}

\subsection{Generalized Deformed Products}

The orthogonal product $P\times Q\subset\R^{d+e}$ of a $d$-polytope $P =
P(A,a) \subset \R^d$ and an $e$-polytope $Q = P(B,b) \subset \R^e$ is given in
inequality description by a block diagonal system:
\[
    \begin{array}{r@{}lcr}
        A x &         & \le & a \\
            &\quad By & \le & b. \\
    \end{array}
\]
We get a \emph{deformed product} (with the combinatorial structure of the
orthogonal product) if we generalize this into a block lower-triangular
system, provided that $Q$ is simple, and that we rescale the right-hand side
of the system suitably. 

\begin{dfn}[Rank $r$ deformed product]
    Let $P = P(A,a) \subset \R^d$ be a $d$-polytope and $Q = P(B,b) \subset
    \R^e$ a simple $e$-polytope, with $A\in\R^{k\times d}$ and $B \in \R^{n
    \times e}$.  Let $C \in \R^{n \times d}$ be an arbitrary matrix of rank
    $r$ and let $M \gg 0$ be large.  The \emph{rank $r$ deformed product} $P
    \bowtie_C Q \subset \R^{d+e}$ of $P$ and $Q$ with respect to $C$ is given
    by 
    \[
    \begin{array}{r@{}lcr}
        A x &     & \le & a \\
        C x &\,+\,By & \le & M\,b \\
    \end{array}, \qquad\textrm{that is,}\qquad
    \left(\begin{array}{cc} A  &\\C &B\end{array}\right)
    \left(\begin{array}{c} x \\y\end{array}\right)\ =\  
    \left(\begin{array}{r} a \\M\,b\end{array}\right).
    \]
\end{dfn}

\begin{prop}\label{prop:is_product}
    Let $P = P(A,a) \subset \R^d$ be a $d$-polytope, $Q = P(B,b) \subset \R^e$
    a simple $e$-polytope, $P \bowtie_C Q$ their deformed product, and $M > 0$
    the parameter involved in its construction. 

    If $M$ is sufficiently large (depending on~$B$, $b$ and $C$), then $P
    \bowtie_C Q$ and $P \times Q$ are combinatorially equivalent.
\end{prop}

Our proposition may also be obtained from the Isomorphism Lemma
\cite[Lemma~2.4]{am98} that was applied by Amenta \& Ziegler to prove the
corresponding statement for (rank 1) deformed products.  However, we use it in
a dual form as given below.  Again, for $I \subseteq [n]$ we write $P_I = P
\cap \{ x : A_I \, x = b_I \}$ for the smallest face $F \subseteq P$ that
satisfies $I \subseteq \eq F$.

\begin{lem}[Isomorphism Lemma; dual formulation]
    Let $P = P(A,a)$ and $Q = P(B,b)$ be two polytopes with $n$ facets and
    $\dim P \ge \dim Q$. If 
    \[
        P_I \text{ is a vertex}\ \ \Longrightarrow\ \ Q_I \text{ is nonempty}
    \]
    for every set $I \subset [n]$ then $P$ and $Q$ are of the same
    combinatorial type.
\end{lem}

\begin{proof}[Proof of Proposition \ref{prop:is_product}] 
    Since $Q$ is a simple polytope, we can find an $M \gg 0$ such that $Q
    \cong P(B, M b - Cv)$ for every $v \in \V P$. In particular, if $u \in \V
    Q$ is a vertex with $I = \eq u$ then $P(B, M b - Cv)_I$ is a vertex.
    Thus, by the dual Isomorphism Lemma, the result follows.
\end{proof}

Proposition \ref{prop:is_product} frees us from a discussion of right hand
sides.  Therefore all deformed products hereafter are understood with a
\emph{suitable} right hand side.

To see that the above definition of rank~$r$ deformed products generalizes the
(rank $1$) deformed products of Amenta \& Ziegler \cite{am98}, we recall their
$\Hr$-description of a deformed product. Let $P= P(A,a) \subset \R^d$ be a
polytope and $\varphi: P \rightarrow \R$ an affine functional with $\varphi(P)
\subseteq [0,1]$.  Let $Q_1,Q_2 \subset \R^e$ be ``normally equivalent''
$e$-polytopes, that is, combinatorially equivalent polytopes with the same
left-hand side matrix, $Q_i = P(B,b_i)$ for $i=1,2$.  Then according to
\cite[Thm.~3.4(iii)]{am98} the exterior representation of $(P,\varphi) \bowtie
(Q_1,Q_2)$ of the \emph{AZ-deformed product} is given by
\[
     (P,\varphi) \bowtie (Q_1,Q_2)\ \  =\  \left\{ (x,y) \in \R^{d+e} : 
     A x \le a,\ B y \le b_1 - (b_1 - b_2) \varphi(x) \right\}
\]

\begin{prop} \label{prop:AZ-prod}
    The AZ deformed product is a rank~$1$ deformed product.
\end{prop}

\begin{proof}
    Let $\varphi(x) = c^\T x + \delta$ be the affine functional.  Let $C =
    (b_1 - b_2) c^\T$ be the matrix of rank at most~$1$ with entries $C_{ij}
    := (b_1 - b_2)_i \cdot c_j$. Further, let $b = b_1 - \delta (b_1 - b_2)$
    and $Q = P(B,b)$. Now, rewriting the inequality system for $(P,\varphi)
    \bowtie (Q_1,Q_2)$ proves the claim.
\end{proof}

\section{Neighborly Cubical Polytopes}\label{sec:ncp}

For $\eps > 0$ the interval $I_\eps = \{ x \in \R : \pm \eps x \le 1 \}$ is a
$1$-dimensional, simple polytope. Its poset of nonempty faces is the poset on
$\{ +, -, 0 \}$ with order relations $+ \prec 0$ and $- \prec 0$. The signs
$\pm$ represent the vertices of the interval with the suggestive notation that
$\pm$ names the vertices given by $\pm\eps x = 1$ while $0$ stands for the
unique (improper) $1$-dimensional face.  An $n$-fold product of intervals
gives a combinatorial $n$-dimensional cube $C_n$ with inequality system
\[
    \begin{array}{r}
        \scriptstyle \pm 1 \\
        \scriptstyle \vdots  \\
        \scriptstyle \pm(n-k) \\
        \scriptstyle \pm(n-k+1) \\
        \scriptstyle \vdots  \\
        \scriptstyle \pm n \\
    \end{array}
    \left(
    \begin{array}{cccccc}
        \pm\eps &         &         &         &         &         \\          
                & \ddots  &         &         &         &         \\
                &         & \pm\eps &         &         &         \\
                &         &         & \pm\eps &         &         \\
                &         &         &         & \ddots  &         \\
                &         &         &         &         & \pm\eps \\
    \end{array}
    \right) x \le 
    \begin{pmatrix} 1 \\ \vdots \\ 1 \\ 1 \\ \vdots \\ 1 \end{pmatrix}.
\]
Every row in the above system represents two inequalities: The $i$-th row
prescribes an upper and a lower bound for the variable $x_i$. Left to the
system are the labels of the rows to which we will refer in the following.

On the level of posets the facial structure is captured by an $n$-fold direct
product of the poset above. The nonempty faces of $C_n$ correspond to the
elements of $\PMZ^n$ with the (component-wise) induced order relation. An
element $\gamma \in \PMZ^n$ represents the unique face $F_\gamma$ with
equality set $\eq F_\gamma = \{ \gamma_i i : i \in [n] \}$ of dimension $\dim
F_\gamma = \#\{ i \in [n] : \gamma_i = 0 \}$. This, in particular, gives the
$f$-vector as $f_i(C_n) = \tbinom{n}{i} 2^{n-i}$.

The cube, as an iterated product of simple $1$-polytopes, lends itself to
deformation beneath the ``diagonal'' that yields, figuratively,  a
\emph{deformed} product of intervals. In the following we construct deformed
cubes that all subscribe to the same deformation scheme. To avoid cumbersome
descriptions, we fix a \emph{template} for a deformed cube.

\begin{dfn}[Deformed Cube Template] \label{con:deformed_cube}
    For $n \ge d \ge 2$, let $G = \{g_1,\dots,g_{d-1}\} \subset \R^{n-d}$ be
    an ordered collection of row vectors and let $\eps > 0$. We denote by
    $C_n(G)$ a deformed cube with lhs matrix
    \begin{equation}\label{eqn:def_cube}
    \newcommand\grow[1]{\multicolumn{4}{c|}{\;\hrulefill\;g_{#1}\;\hrulefill\;}}
    \newcommand\vrow{\multicolumn{4}{c|}{\vdots}}
    A(G) = (\ovA,\oovA) = 
    \left(
    \begin{array}{cccc|cccc}
     \pm\eps    &         &         &         &         &         &        &  \\
     1          & \pm\eps &         &         &         &         &        &  \\
                & 1       & \ddots  &         &         &         &        &  \\
                &         & \ddots  & \pm\eps &         &         &        &  \\
                &         &         & 1       & \pm\eps &         &        &  \\
     \grow{1}                                 &         & \pm\eps &        &  \\
     \vrow                                    &         &         & \ddots &  \\
     \grow{d-1}                               &         &         &&\pm\eps   \\
    \end{array}
    \right).
    \end{equation}
\end{dfn}

Proposition \ref{prop:is_product} assures of a suitable right hand side such
that $C_n(G)$ is a combinatorial $n$-cube. Up to this point, we required
$\eps$ to be nothing but positive; this will be subject to change, soon. 

The polytope we are striving for is the image of $C_n(G)$ under projection.
Recall that our projections will be onto the last $d$ coordinates for which
the vertical bar in (\ref{eqn:def_cube}) is a reminder.

We now come to the first main result of this section.

\begin{thm}[Joswig \& Ziegler {\cite[Theorem 17]{jos00}}]\label{thm:ncp}
    For every $2 \le d \le n$ there is a cubical $d$-polytope whose $( \lfloor
    \frac{d}{2} \rfloor - 1)$-skeleton is isomorphic to that of an $n$-cube.
\end{thm}

\begin{proof}
    The claim will be established by choosing the \emph{right} deformation and
    verifying that all the necessary faces are strictly preserved under
    projection.

    Let $Q$ be a neighborly $(d-2)$-polytope with $n-1$ vertices in general
    position. In particular, $Q$ has the property that every subset of at most
    $\lfloor\frac{d-2}{2}\rfloor = \lfloor \frac{d}{2} \rfloor - 1$ vertices
    forms a face of $Q$.  For an arbitrary but fixed ordering of the vertices,
    let $G \in \R^{(n-1)\times (n-d)}$ be a Gale transform of $Q$.

    As the vertices of $Q$ are in general position, we can choose a Gale
    transform of the form $G = \left( \begin{smallmatrix} I_{n-d} \\ \ovG
    \end{smallmatrix} \right)$, where $\ovG = \{g_1,\dots,g_{d-1}\} \subset
    \R^{n-d}$ is an ordered collection of row vectors. Let $C = C_n(\ovG)$ be
    the deformed cube given by the template (\ref{eqn:def_cube}) with respect
    to $\ovG$.

    We claim that the projection of $C$ to the last $d$ coordinates yields the
    result. For this we prove that all faces of dimension up to $k = \lfloor
    \frac{d}{2} \rfloor -1$ survive the projection. In order to do so, we
    propose the following strategy: We will show that for an arbitrary vertex
    $v$ of $C$ the incident faces of dimension $\le k$ are retained.  

    Consider $\ovA_{\eq v}$, the first $n-d$ columns of the inequalities of
    (\ref{eqn:def_cube}) which are tight at $v$. The matrix is of the form
    \begin{equation}\label{eqn:Aeq}
       \ovA_v := 
       \ovA_{\eq v } = 
        \left(
        \begin{array}{cccc}
            \s_1\eps &          &        &              \\
            \hline
              1      & \s_2\eps &        &              \\
                     & 1        & \ddots &              \\
                     &          & \ddots & \s_{n-d}\eps \\
                     &          &        & 1            \\
            \multicolumn{4}{c}{\;\hrulefill\;g_1\;\hrulefill\;}\\
            \multicolumn{4}{c}{\vdots}     \\
            \multicolumn{4}{c}{\;\hrulefill\;g_{d-1}\;\hrulefill\;}   \\
        \end{array}
        \right) \in \R^{n \times (n-d)}
    \end{equation}
    with $\s_1,\dots,\s_{n-d} \in \PM$.

    Since the vertices of $Q$ are in general position, by Proposition
    \ref{prop:simplicial_gale}, $G$ is a configuration of vectors in general
    position with respect to linear hyperplanes. Thus, for $\eps > 0$
    sufficiently small, $\ovA_v$ take away the first row is still the Gale
    transform of a polytope combinatorially equivalent to $Q$. By Gale
    duality, this in particular means that discarding up to $\lfloor
    \frac{d-2}{2} \rfloor = k$ rows from $\ovA_v$ leaves the remaining ones
    positively spanning.
    
    Now, let $F \subset C$ be a face of dimension $\ell \le k$ with $v \in F$.
    By the Projection Lemma \ref{lem:projection}, $F$ is strictly preserved by
    the projection iff the rows of $\ovA_I$ for $I = \eq F$ are positively
    spanning. Since $C$ is simple, $\ovA_I$ is an $n-\ell$ rowed submatrix of
    $\ovA_{\eq v}$, that is, at most $k$ rows have been discarded from
    $\ovA_v$.

    Choosing $\eps$ sufficiently small also has the effect that the rows of
    $\ovA_v$ are in general position with respect to linear hyperplanes.
    Thus, Corollary \ref{cor:all_strictly} vouches for the fact that all faces
    of $\pi_d(C_n(\ovG))$ arise from the projection of~$C_n(\ovG)$.
\end{proof}

The polytope $\pi_d(C_n(\ovG))$ constructed in the course of the proof depends
on the choice of a neighborly $(d-2)$-polytope $Q$ with $n-1$ vertices in
general position, equipped with an ordering of its vertices.  In particular,
the order of the vertices is needed to determine $\ovG$ and thus $C_n(\ovG)$.
Nevertheless, by abuse of notation we will write $C_n(Q)$ for the deformed
cube $C_n(\ovG)$. We will see in the next section that, in fact, the
combinatorics of $\pi_d(C_n(Q))$ is determined  by the choice of $Q$ and the
vertex order. In Section \ref{sec:ncp_cyc}, we show that the polytopes
constructed in \cite{jos00} correspond to the case were $Q$ is a cyclic
polytope with the standard vertex ordering.  For now, we baptize the polytope
that we have constructed.

\begin{dfn}
    For parameters $n \ge d \ge 2$ and a neighborly $(d-2)$-polytope $Q$ on
    $n-1$ ordered vertices in general position, we denote the \emph{neighborly
    cubical polytope} $\pi_d(C_n(Q))$ by $\ncp_{n,d}(Q)$.
\end{dfn}

Let us briefly comment on the extremal  choices of $d$.  For $d = n$, the
polytope $\ncp_{n,n}(Q)$ is combinatorially isomorphic to an $n$-cube. The
neighborly polytope $Q$ is then an $(n-2)$-polytope with $n-1$ vertices, a
simplex.  For $d = 2$, the polytope $\ncp_{n,2}(Q)$ is a $2^n$-gon and
$C_n(Q)$ is, in fact, a realization of a \emph{Goldfarb cube} \cite{gold82}.
What might strike the reader as strange is that the neighborly polytope in
question is a $0$-dimensional polytope with $n-1$ vertices. The Gale transform
of such a polytope is given by the vertices of a $(n-2)$-simplex with vertices
$\{ e_1, e_2, \dots, e_{n-2}, -\1 \}$. 

The proof can be adapted to yield a \emph{$k$-neighborly cubical polytope},
that is, a polytope having its $k$-skeleton isomorphic to that of an $n$-cube.
By \cite[Corollary 5]{jos00}, the neighborliness is bounded by $k \le \lfloor
\frac{d}{2} \rfloor - 1$. In our construction this fact is reflected as
follows. The polytope $\ncp_{n,d}(Q)$ is $k$-neighborly cubical iff $Q$ is
$k$-neighborly.  By \cite[Exercise 0.10]{zie95}, neighborliness for
$(d-2)$-polytopes is bounded by $\lfloor \frac{d - 2}{2} \rfloor$. 

\subsection{Combinatorial description of the neighborly cubical polytopes}
\label{sec:ncp_lattice}

We describe the face lattice of $\ncp_{n,d}(Q)$ in terms of
\emph{lexicographic triangulations} of $Q$. We start by giving the necessary
background on regular subdivisions with an emphasis on lexicographic
triangulations in terms of Gale transforms. Our main sources are the paper by
Lee \cite{lee91} and the (upcoming) book by De Loera \emph{et al.} \cite{LRS}.

Let $Q$ be a simplicial $D = d-2$ dimensional simplicial polytope on $N = n -
1$ ordered vertices. We further assume that the vertices of $Q$ are in general
position, i.e.\ all vertex induced subpolytopes are simplicial as well. Let
the rows of $V \in \R^{N \times D}$ be the vertices of~$Q$ in some ordering,
and let $\w = (\w_1,\dots,\w_N)^\T \in \R^N$ be a set of \emph{heights}.
Denote by $V^\w = ( \w, V ) \in \R^{N \times (D+1)}$ the ordered set of
\emph{lifted} vertices $(\w_i,v_i)$ for $i = 1, \dots, N$.  Let $a =
(\w_0,v_0) \in \R^{D+1}$ be arbitrary with $\w_0 \gg \max_i |\w_i|$ and
consider the polytope $Q^\w = \conv ( V^\w \cup a )$. If $\w_0$ is
sufficiently large, then the vertex figure of $a$ in $Q^\w$ is isomorphic to
$Q$ and the closed star of $a$ in $\partial Q^\w$ is isomorphic to that of the
apex of a pyramid over $Q$.  The anti-star (or deletion) of $a$ in the
boundary of $Q^\w$, i.e.\ the faces of $Q^\w$ not containing $a$, constitute a
pure $D$-dimensional polytopal complex $\Gamma_\w$, the \emph{$\w$-induced}
(or \emph{$\w$-coherent}) \emph{subdivision}. The name ``subdivision'' stems
from the fact that the underlying set $\|\Gamma_\w\|$ is piecewise-linear
homeomorphic to $Q$ via the projection onto the last $D$ coordinates.  The
inclusion maximal polytopes in $\Gamma_\w$ are called \emph{cells}.
$\Gamma_\w$ is called a \emph{triangulation} if every cell is a $D$-simplex.
Altering the heights $\w^\prime_i = \w_i + \ell(v_i)$ along an affine
functional $\ell: Q \rightarrow \R$ leaves the induced subdivision unchanged.
We call a set of heights \emph{normalized} if its support is minimal in the
corresponding equivalence class.

\begin{prop}\label{prop:def_pyr}
    Let $\w^\T = (\w_1,\dots,\w_{N-D-1}, 0, \dots, 0) \in \R^N$ be a
    normalized set of heights and let $G = \left( \begin{smallmatrix}
    \id{N-D-1} \\ \ovG \end{smallmatrix} \right) \in \R^{N \times (N-D-1)}$.
    For $\eps > 0$ sufficiently small, the matrix
    \[
        G_\w =
            \begin{pmatrix}
                -\eps \ovw \\ \;\;\,G 
            \end{pmatrix}
        \in \R^{(N+1) \times (N - D-1)}
    \]
    with $\ovw = (\w_1,\dots,\w_{n-d-1})$ is a Gale transform of a polytope
    combinatorially equivalent to $Q^\w$.
\end{prop}
\begin{proof}
    It is easily verified that the columns of 
    \[
        \left(
        \begin{array}{ccc}
            1 & 1 & \0 \\
            \1 + \eps \w & \eps \w & V \\
        \end{array}
        \right) \in \R^{(N+1) \times (D+2)}
    \]
    form a basis for the orthogonal complement of the column span of $G_\w$.
    For $\eps$ sufficiently small, the first column is strictly positive and
    dehomogenizing with respect to this column yields the desired polytope.
\end{proof}
In particular, $G_\w$ encodes the combinatorial structure of $Q$ as well as
that of the $\w$-induced regular subdivision.

Consider the two induced regular subdivisions of~$Q$ obtained by lifting the
vertex $v_1$ to height $\w_1 = \pm h$ with $h > 0$ and fixing all the
remaining heights to $0$. In both cases the lifted polytope is a pyramid over
the polytope $Q^\prime = \conv( V \setminus v_1)$.  For $\w_1 = -h$ the
subdivision is said to be obtained by \emph{pulling} $v_1$ and its cells are
pyramids over the remote facets of $Q^\prime$, that is, the facets common to
both $Q$ and $Q^\prime$.  This subdivision is, in fact, a triangulation since
its cells are pyramids over $(D-1)$-simplices.  The other subdivision ($\w =
+h$) is said to be obtained by \emph{pushing} $v_1$ and its cells are pyramids
over the newly created facets of $Q^\prime$, which are again simplices, plus
one (possibly non-simplex) cell that is $Q^\prime$.  

The ordering of the vertices of $Q$ gives rise to a chain of (sub-)polytopes
$Q = Q_0 \supset Q_1 \supset \cdots \supset Q_{N-D-1} = \Delta_D$ with $Q_i =
\conv \{ v_{i+1},\dots,v_N \}$ simplicial $D$-polytopes. Let $1 \le k \le N -
D - 1$, then the \emph{$k$-th lexicographic triangulation} $\Lex{k} Q$ of $Q$
in the given vertex order is the triangulation obtained by pushing the first
$k-1$ vertices in the given order and then pulling the $k$-th vertex. That is
to say, pushing $v_1$ creates a subdivision of $Q = Q_0$ that has $Q_1$ as its
only non-simplex cell.  Subsequently, the cell $Q_1$ gets replaced by a
pushing subdivision of $Q_1$ with respect to $v_2$, and so on. Finally,
pulling $v_{k+1}$ in $Q_k$ completes the triangulation. The following lemma
asserts that the above procedure yields a regular subdivision by giving a
description in the spirit of Proposition \ref{prop:def_pyr}.

\begin{lem}[{\cite[Example 2]{lee91} \cite{dipl-Sanyal}}]\label{lem:lex_gale}
    Let $\eps > 0$ and $\w = (\w_1,\w_2,\dots,\w_{N-D-1},0,\dots,0)\in \R^N$
    be a set of normalized heights satisfying $|\w_{i+1}| \le \eps |\w_i|$ for
    all $1 \le i \le N-D-2$. If $\eps > 0$ is sufficiently small, then $G_\w$
    is a Gale transform encoding $\Lex{k}Q$ for
    \[
        k \ \ =\ \ \min \{i : \w_i < 0 \} \cup \{ n-d-1 \}.
    \]
\end{lem}

\begin{dfn}
    We call the polytope $\L_k(Q) = \Qw^\w$ corresponding to $G_\w$ the
    \emph{$k$-th lexicographic pyramid} of $Q$. 
\end{dfn}

According to the remarks following Proposition \ref{prop:def_pyr}, $\L_k(Q)$
carries both the combinatorics of $Q$ as well as that of $\Lex{k}Q$. So every
facet of $\L_k(Q)$ is either a pyramid over a facet of $Q$ or a cell of
$\Lex{k}Q$.

We are now in a position to determine the combinatorics of $\ncp_{n,d}(Q)$.
To be more precise, we determine the \emph{local} combinatorial structure,
i.e.\ for any given vertex we describe the set of facets that contain it. The
construction of a neighborly cubical polytope depended on an ordering of the
vertices of $Q$, which we fix for the following theorem.

\begin{thm}\label{thm:comb_ncp}
    Let $C = C_n(Q)$ be the deformed cube with respect to $Q$. Further, let $v
    \in C$ be an arbitrary vertex with $\eq v$, given by  $\s \in \PM^n$.
    Then the vertex figure of $\pi_d(v)$ in $\ncp_{n,d}(Q)$ is isomorphic
    to~$\L_p(Q)$ for
    \[
        p \ \ =\ \ \min \{ i \in [n] : \s_i = + \} \cup \{n-d-1\}.
    \]
    In particular, the $(d-1)$-faces of $C$ containing $v$ that are preserved
    by projection are in one-to-one correspondence to the facets of $\L_p(Q)$.
\end{thm}

\begin{proof}
    After a suitable base transformation of (\ref{eqn:Aeq}) by means of column
    operations, the first $n-d$ columns of $A_{\eq v}$ can be assumed to be of
    the form
    \[
        \newcommand\g{\tilde{g}}
        \left(
        \begin{array}{cccc}
                -\w_1 &     -\w_2 & \cdots & -\w_{n-d}     \\
            \hline
              1      &          &        &              \\
                     & 1        &        &              \\
                     &          & \ddots &              \\
                     &          &        & 1            \\
            \multicolumn{4}{c}{\;\hrulefill\;\g_1\;\hrulefill\;}\\
            \multicolumn{4}{c}{\vdots}     \\
            \multicolumn{4}{c}{\;\hrulefill\;\g_{d-1}\;\hrulefill\;}   \\
        \end{array}
        \right) 
    \]
    with
    \[
        \w_i \ \ =\ \ (-1)^i\eps^i \prod_{j=1}^i\s_j
    \]
    By Lemma \ref{lem:lex_gale}, this is a Gale transform of $\L_k(Q)$ with $k
    = p$. 

    Any generic projection of polytopes $\pi : P \rightarrow P^\prime =
    \pi(P)$ induces a (contravariant) order and rank preserving map $\pi^\#:
    \FL \partial P^\prime \hookrightarrow \FL \partial P$.

    The face poset of $\partial\ncp_{n,d}(Q)/u$, the boundary complex of the
    vertex figure of $u = \pi_d(v)$ in $\ncp_{n,d}(Q)$, is isomorphic to
    $\pi^\#(\FL \partial\ncp_{n,d}(Q)_{\ge u})$, the image of the principal
    filter of $u$. By the Projection Lemma, the image coincides with the
    embedding of $\L_p(Q)$ into the vertex figure $\FL \partial C_n(Q)_{\ge
    v}$.
\end{proof}

Theorem \ref{thm:comb_ncp} implies that the quotient $\ncp_{n,d}(Q) / e$  with
respect to certain edges is isomorphic to~$Q$. This observation implies the
following result.

\begin{cor}
     Non-isomorphic neighborly $(d-2)$-polytopes $Q$ and $Q^\prime$ yield
     non-isomorphic neighborly cubical polytopes $\ncp_{n,d}(Q)$ and
     $\ncp_{n,d}(Q^\prime)$. Moreover, there are at least as many different
     combinatorial types of $d$-dimensional neighborly cubical polytopes as
     there are neighborly simplicial $(d-2)$-polytopes on $n-1$ vertices.
\end{cor}

The number of combinatorial types of neighborly simplicial polytopes is huge,
according to Shemer~\cite{shemer}.

\subsection{Neighborly cubical polytopes from cyclic polytopes}
\label{sec:ncp_cyc}

In this section we (re-)construct the neighborly cubical polytopes of Joswig
\& Ziegler \cite{jos00}. This specializes the discussion in the previous
section to the case of  $Q$ a \emph{cyclic polytope} in the standard vertex
ordering. By a thorough analysis of the lexicographic triangulations of cyclic
polytopes we recover the ``cubical Gale's evenness criterion'' of
\cite{jos00}. For a treatment  of cyclic polytopes and their triangulations
beyond our needs we refer the reader to \cite{LRS} and \cite{zie95}.

The degree $D$ \emph{moment curve} is given by $t \mapsto \gamma(t) =
(t,t^2,\dots,t^D) \in \R^D$. For given pairwise distinct values $t_1, t_2,
\dots, t_N \in \R$ with $N \ge D+1$ the convex hull of the corresponding
points on the moment curve $\cyc_D(t_1,\dots,t_N) = \conv \{ \gamma(t_i) : i
\in [N] \}$ is a convex $D$-dimensional polytope. A fundamental consequence of
the theorem below is that the combinatorial type of $\cyc_D(t_1,\dots,t_N)$ is
independent of the actual values $t_i$. Therefore, we work with $\cyc_D(N) :=
\cyc_d(1,2,\dots,N)$, the $D$-dimensional \emph{cyclic polytope} on $N$
vertices in standard order. For the sake of notational convenience later on,
we describe its faces in terms of characteristic vectors of cofaces: A vector
$\a \in \ZO^N$ names a coface of $\cyc_D(N)$ iff $\conv \{ \gamma(i) : \a_i =
0 \}$ is a face of $\cyc_D(N)$. We also extend the notion of ``co-'' to
subdivisions and, therefore, speak freely about \emph{cocells}.

Let $\a \in \ZO^N$ such that $\#\{ j < i : \a_j = 0 \}$ has the same parity
for every $i \in [N]$ with $\a_i = 1$. Then $\a$ is called \emph{even} or
\emph{odd} according to this parity.

\begin{thm}[Gale's Evenness Criterion {\cite[Sect.~4.7]{grue03}
    \cite[Thm.~0.7]{zie95} \cite[Thm.~6.2.6]{LRS}}]
    A vector $\a \in \ZO^N$ names a cofacet of $\cyc_D(N)$ if and only if $\a$
    has exactly $D$ zero entries and is either even or odd.
\end{thm}

As a byproduct we get that cyclic polytopes are
\begin{compactitem}[$\bullet$~]
    \item simplicial, since all facets have exactly $D$ vertices,
    \item in general position, since every subpolytope is again cyclic, and
    \item neighborly, since every $\a \in \ZO^N$ with $\le \lfloor \frac{D}{2}
        \rfloor$ zeros can be made to meet the above conditions by changing
        entries $1 \rightarrow 0$.
\end{compactitem}

From a geometric point of view, the odd and even (co)facets correspond to the
upper and lower facets of $\cyc_D(N)$ with respect to the last coordinate.
This dichotomy among the facets allows for an explicit characterization of the
(simplicial) cells of a pushing/pulling subdivision of $\cyc_D(N)$ with
respect to the first vertex. Moreover, since every vertex induced subpolytope
of $\cyc_D(N)$ is again cyclic and from this we will derive a complete
description of the lexicographic triangulations of cyclic polytopes with
vertices in standard order.

To prepare for the precise statement, let $Q = \cyc_D(N) = \conv \{ v_i =
\gamma_d(i) : i \in [N] \}$ and $Q^\prime = \conv \{ v_2,\dots,v_N \} \cong
\cyc_D(N-1)$ the subpolytope on all vertices except the first.  Let $\Gamma$
be the subdivision of $Q$ obtained by pulling or pushing $v_1$.  Any cell in
$\Gamma$ that contains $v_1$ is a $D$-simplex and, therefore, let $\a \in
\ZO^N$ be a cocell with $D+1$ zero entries and $\a_1 = 0$.  Indeed, any such
cell is a pyramid over a facet of $Q^\prime$ and thus $\a$ is of the form $\a
= (0,\a^\prime)$ and $\a^\prime$ adheres to the Gale's evenness criterion. The
cocell $\a$ is part of a pushing or a pulling subdivision of $Q$ if and only
if $\a$ is or is not a cofacet of $Q$.  Clearly, the first gap in $\a$ is even
and, hence, the parity of the gaps of $\a^\prime$ concludes the
characterization.

\begin{lem}\label{lem:cyc_lex}
    Let $Q = \cyc_D(N)$ be a cyclic polytope and let $\L_k(Q)$ be a
    lexicographic pyramid of $Q$.  Let $\a \in \ZO^{N+1}$ with $D+1$ zero
    entries and let $p = \min\{ i : \a_i = 0 \}$. Thus $\a$ is of the form 
\vskip-5mm
    \[
        \a = (\underbrace{1,1,\dots,1}_{p-1},0,\a^\prime).  
    \]
    Then $\a$ is a cofacet of $\L_k(Q)$ if and only if one of the following
    conditions is satisfied:
    \begin{compactenum}[\rm i)]
        \item $1 = p$ and $\a^\prime$ is a cofacet of
                  $\cyc_d(n)$.
        \item $1 < p < k $ and $\a^\prime$ is even.
        \item $p = k $ and $\a^\prime$ is odd.
    \end{compactenum}
\end{lem}
\begin{proof}
    Every facet containing the $0$-th vertex is a pyramid over a facet of $Q$
    and every incident facet is of the form $\a = (0,\a^\prime)$ with
    $\a^\prime$ a cofacet of $Q$.

    If $2 \le p < k$ then $\a$ names a cocell of the pushing subdivision of
    $Q_{p-1} = \conv\{ v_p, \dots, v_N \}$ with respect to $v_p$ and
    containing $v_p$. This, however, is the case if and only if $\a^\prime$ is
    an even cofacet of~$Q_{p}$.  The case $p = k$ follows from similar
    considerations.
\end{proof}

Setting $N = n-1$ and $D = d-2$ and combining the above description with
Theorem \ref{thm:comb_ncp}, we obtain the following result of Joswig \&
Ziegler.

\begin{thm}[Cubical Gale's Evenness Condition \cite{jos00}]
    Let $F$ be a $(d-1)$-face of the deformed cube $C = C_n(Q)$ with $Q
    :=\cyc_{d-2}(n-1)$. Let $\eq F$ be given by $\a \in \PMZ^n$ and let $p \ge
    1$ be the smallest index such that $\a_{p} = 0$.  The face $F$ projects to
    a facet of $\ncp_{n,d}(Q)$ if and only if $\a$ is of the form
    \[
        \a = (\underbrace{-,-,\cdots,-}_{p-2},\s,0,\a^\prime)
    \]
    with $|\a^\prime| = (|\a^\prime_{p+1}|,\dots, |\a^\prime_n|) \in
    \ZO^{n-p}$ satisfies the ordinary Gale's evenness condition and for $p >
    1$ one of the following conditions holds:
    \begin{compactenum}[\rm i)]
        \item $\s = -$ and $|\a^\prime|$ is even, or
        \item $\s = +$ and $|\a^\prime|$ is odd.
    \end{compactenum}
\end{thm}

\begin{proof}
    Let $v \in F \subset C$ be a vertex with equality set $\b = \eq v$ and
    such that $\b_{p} = +$. By Theorem~\ref{thm:comb_ncp}, the vertex figure
    of $\pi_d(v)$ in $\ncp_{n,d}(Q)$ is isomorphic to $\L_k(Q)$, with $k \in
    \{p-1,p\}$.
    
    Thus $F$ projects to  a facet of $\ncp_{n,d}(Q)$ if and only if $|\a|$ is
    a cofacet of $\L_k(Q)$. The result now follows from
    Lemma~\ref{lem:cyc_lex} by noting that $k = p-1$ iff $\s = -$.
\end{proof}

\section{Deformed Products of Polygons}
\label{sec:pdpp}

The \emph{projected deformed products of polygons} (PDPPs) are $4$-dimensional
polytopes. They were constructed in \cite{zie04} because of their extremal
$f$-vectors: For these polytopes the \emph{fatness} parameter
$\Phi(P):=\frac{f_1+f_2-20}{f_0+f_3-10}$ is large, getting arbitrarily close
to~$9$.  This parameter, introduced in~\cite{zie02}, is crucial for the
$f$-vector theory of $4$-polytopes.  In \cite{zie04} the $f$-vectors of the
PDPPs were computed without having a combinatorial characterization of the
polytopes in reach.

However, the PDPPs are yet another instance of projections of deformed
products, so the theory developed here gives us a firm grip on their
properties. In the following we generalize the construction to higher
dimensions and analyze its combinatorial structure using the tools developed
in this paper. In particular, a description of the facets of the PDPPs appears
for the first time.

To begin with, the following is a generalization of Theorem \ref{thm:ncp}.
\begin{thm}\label{thm:dpp}
    Let $\gon \ge 4$ be even. For every $2 \le d \le 2r$ there is a
    $d$-polytope whose $(\lfloor \frac{d}{2}\rfloor-1)$-skeleton is
    combinatorially isomorphic to that of an $r$-fold product of $\gon$-gons.
\end{thm}

Let us remark that the proofs of the results in this section can be adapted to
yield the generalizations for products of even polygons with varying numbers
of vertices in each factor. However, the generalized results require more
technical and notational overhead. Therefore, we trade generality in for
clarity and only give the uniform versions of the results.

For $\gon = 4$ the $r$-fold product of quadrilaterals is actually a cube of
dimension $n = 2r$ and thus  $\ncp_{n,d}(Q)$ satisfies the claims made. In the
inequality description the quadrilaterals can be seen by pairing up the
intervals indicated by the framed submatrices below:

\begin{equation*} \label{eqn:ncp-quad}
    \newcommand{\lB}[1]{\multicolumn{1}{|@{}c@{}}{\,#1\,}}
    \newcommand{\rB}[1]{\multicolumn{1}{@{}c@{}|}{\,#1\,}}
    \newcommand{\E}{\pm\eps}
    \left(\,
    \begin{array}{ccccccc|ccccc}
\cline{1-2}
    \lB{\E}& \rB{ }  &         &         &         &         &         &         &         &        &         &        \\
    \lB{1} & \rB{\E} &         &         &         &         &         &         &         &        &         &        \\
\cline{1-2} \cline{3-4}
           &   1     & \lB{\E} & \rB{ }  &         &         &         &         &         &        &         &        \\
           &         & \lB{1}  & \rB{\E} &         &         &         &         &         &        &         &        \\
\cline{3-4}
           &         &         &         & \ddots  &         &         &         &         &        &         &        \\
\cline{6-7} 
           &         &         &         &         & \lB{\E} & \rB{ }  &         &         &        &         &        \\
           &         &         &         &         & \lB{1}  & \rB{\E} &         &         &        &         &        \\
\cline{6-7} \cline{8-9}
           &         &         &         &         &         &  1      & \lB{\E} & \rB{ }  &        &         &        \\
\multicolumn{7}{c|}{\hrulefill\;g_{1}\; \hrulefill}                  & \lB{ }  & \rB{\E} &        &         &          \\
\cline{8-9}
\multicolumn{7}{c|}{           \vdots           }                     &         &         & \ddots &         &         \\
\cline{11-12}
\multicolumn{7}{c|}{\hrulefill\;g_{d-2}\;\hrulefill}                  &         &         &        & \lB{\E} & \rB{ } \\
\multicolumn{7}{c|}{\hrulefill\;g_{d-1}\;\hrulefill}                  &         &         &        & \lB{ }  & \rB{\E}\\
\cline{11-12}
    \end{array}        
    \,\right).
\end{equation*}

We wish to build on this special case and therefore consider the normals of
such a quad:
\[
    \psfrag{e1}[lc]{\scriptsize{$(+\eps,0)$}}
    \psfrag{e2}[rc]{\scriptsize{$(-\eps,0)$}}
    \psfrag{e3}[lb]{\scriptsize{$(1,+\eps)$}}
    \psfrag{e4}[lb]{\scriptsize{$(1,-\eps)$}}
    \includegraphics[width=2cm]{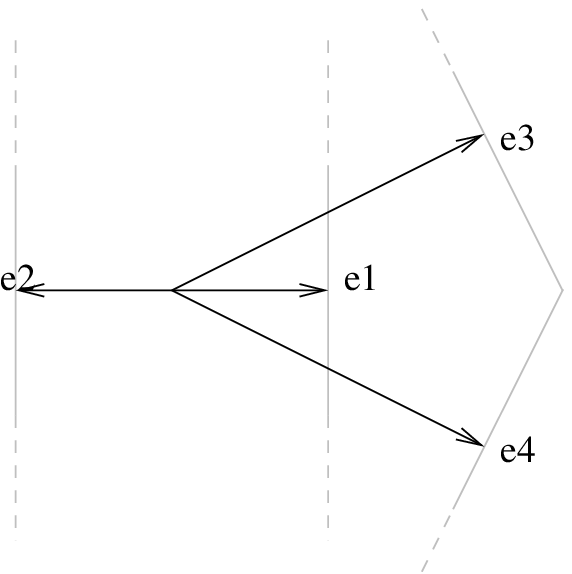}
\]

The polygons we are heading for arise as generalizations of the above quad. For
$\gon \ge 4$ even, consider the vectors
\[
\begin{array}{r@{\,=(}rr@{}l}   
        a_0  &  -1,& 0 )\\[4pt]
        a_i  &   1,& \eps \tfrac{\gon-2i}{\gon-2} ) & \quad\text{
        for } i = 1,\dots,\gon-1\\
    \end{array}   
\]
as shown below.  For suitable $b_0, b_1, \dots, b_{\gon-1} > 0$, 
\[
    a_i^\T x \le b_i \text{ for } i=0, \dots, \gon-1
\]
describes a convex $\gon$-gon in the plane:
\[
    \psfrag{L}[cb]{}
    \psfrag{Q}[cb]{}
    \psfrag{D}{}
    \psfrag{a0}[rc]{\scriptsize{$(-1,0) = a_0$}}
    \psfrag{a1}[lb]{\scriptsize{$a_1=(1,+\eps)$}}
    \psfrag{a2}[lc]{\scriptsize{$a_2$}}
    \psfrag{ak-1}[lc]{}
    \psfrag{ak}[lb]{\scriptsize{$a_{\gon/2}=(1,0)$}}
    \psfrag{ak+1}[lc]{}
    \psfrag{an-2}[lc]{\scriptsize{$a_{\gon-2}$}}
    \psfrag{an-1}[lc]{\scriptsize{$a_{\gon-1}=(1,-\eps)$}}
    \includegraphics[width=4.5cm]{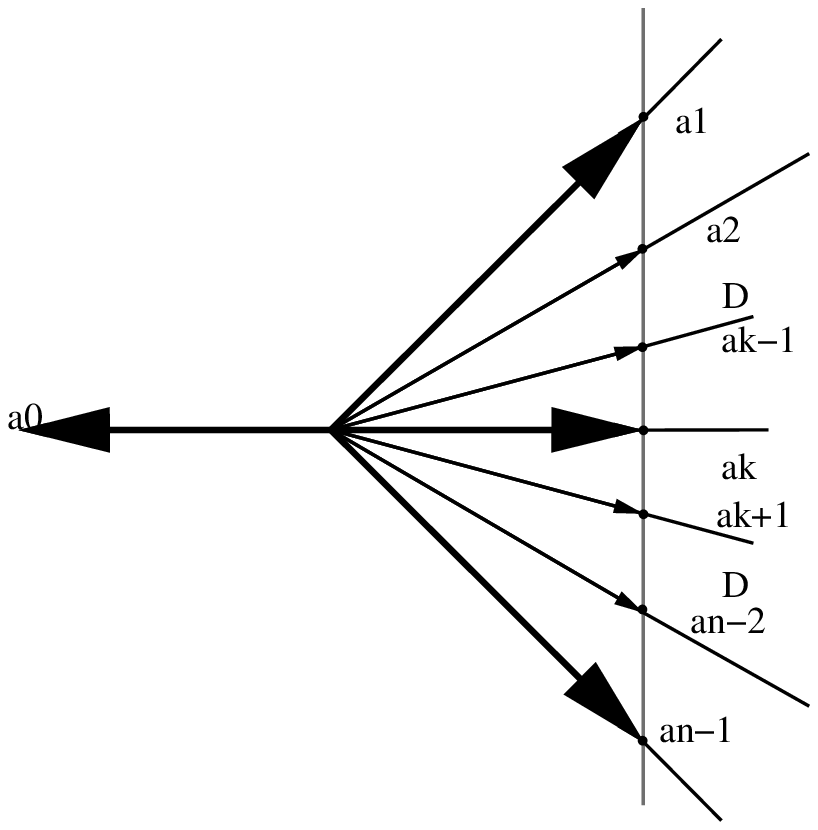}
\]
For the finishing touch, we scale every even-indexed inequality by $\eps$,
\[
    \psfrag{D}{\scriptsize{$\vdots$}}
    \psfrag{a0}[rb]{\scriptsize{$(-\eps,0) = \eps a_0$}}
    \psfrag{a1}[lb]{\scriptsize{$a_1$}}
    \psfrag{a2}[lc]{\scriptsize{$\eps a_2$}}
    \psfrag{ak-1}[lc]{\scriptsize{$a_{\gon/2-1}$}}
    \psfrag{ak}[lc]{\scriptsize{$\eps a_{\gon/2}$}}
    \psfrag{ak+1}[lc]{\scriptsize{$a_{\gon/2+1}$}}
    \psfrag{an-2}[lc]{\scriptsize{$\eps a_{\gon-2}$}}
    \psfrag{an-1}[lc]{\scriptsize{$a_{\gon-1}$}}
    \includegraphics[width=4.5cm]{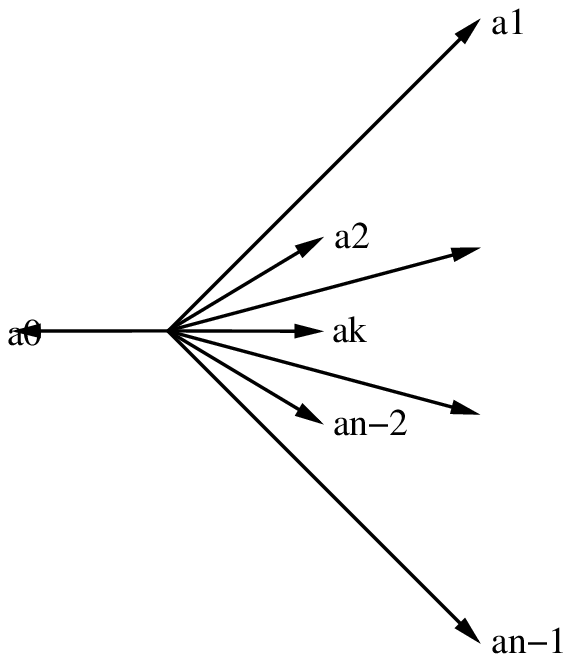}
\]
We arrange the scaled normals  and right
hand sides into a matrix and vector respectively:
\[
    A = \left(
    \begin{array}{r@{}l}
        \eps & a_0 \\
             & a_1 \\
        \eps & a_2 \\
        \multicolumn{2}{c}{\vdots}\\
             & a_{\gon-1} \\
    \end{array}
    \right)
    \ \ 
    \text{ and } 
    \ \ 
    b = \left(
    \begin{array}{r@{}l}
        \eps & b_0 \\
             & b_1 \\
        \eps & b_2 \\
        \multicolumn{2}{c}{\vdots}\\
             & b_{\gon-1} \\
    \end{array}
    \right).
\]

Using these \emph{special} polygons we set up a template for a deformed
product of polygons (DPP).

\begin{dfn}[DPP Template]
    For $\gon \ge 4$ even and $2r \ge d \ge 2$, let $G = \{g_1,\dots,g_{d-1}\}
    \subset \R^{2r-d}$ be an ordered collection of row vectors. We denote by
    $\dpp{2r}{G}$ the deformed product of polygons with lhs inequality system 
    \begin{equation} \label{eqn:dpp}
    \newcommand{\Bx}[2]{
        \ifthenelse{#1>0}
        {
            \multicolumn{2}{|c|}{\raisebox{0.7ex}[0ex]{\Large$A_n$}}
        }
        {
            \multicolumn{2}{|c}{\raisebox{0.7ex}[0ex]{\Large$A_n$}}
        }
    }
    \newcommand{\Bl}{\multicolumn{2}{|c|}{ }}
    \newcommand{\Bbl}{\multicolumn{2}{|c}{ }}
    \newcommand{\A}{\multicolumn{2}{|c|}{\raisebox{0.7ex}[0ex]{\Large$A$}}}
    \newcommand{\Aa}{\multicolumn{2}{|c}{\raisebox{0.7ex}[0ex]{\Large$A$}}}
    \newcommand{\vl}{\vline\,\vline\;}
    \newcommand{\bl}{\multicolumn{1}{|c|}{ }}
    \newcommand{\lB}[1]{\multicolumn{1}{|@{}c@{}}{\,#1\,}}
    \newcommand{\rB}[1]{\multicolumn{1}{@{}c@{}|}{\,#1\,}}
    \newcommand{\rBB}[1]{\multicolumn{1}{@{}c@{}||}{\,#1\,}}
    \left(\,
    \begin{array}{ccccccc|ccccc}
\cline{1-2}
        \Bl          &         &         &         &         &         &         &         &        &         &       \\
        \A           &         &         &         &         &         &         &         &        &         &       \\
\cline{1-2} \cline{3-4}
     \lB{ }&   1     & \Bl               &         &         &         &         &         &        &         &       \\
     \lB{ }&         & \A                &         &         &         &         &         &        &         &       \\
\cline{1-2} \cline{3-4}
           &         &         &         & \ddots  &         &         &         &         &        &         &       \\
\cline{6-7} 
           &         &         &         &         & \Bl               &         &         &        &         &       \\
           &         &         &         &         & \A                &         &         &        &         &       \\
\cline{6-7} \cline{8-9}
\cline{1-9}
    \lB{ } &         &         &         &         &         &  1      & \Bl               &        &         &       \\
\multicolumn{7}{|c|}{\hrulefill\;g_{1}\; \hrulefill}                   & \A                &        &         &       \\
\cline{1-9}
\multicolumn{7}{c|}{           \vdots           }                      &         &         & \ddots &         &       \\
\cline{11-12}
\cline{1-7}
\multicolumn{7}{|c|}{\hrulefill\;g_{d-2}\;\hrulefill}                  &         &         &        & \Bl            \\
\multicolumn{7}{|c|}{\hrulefill\;g_{d-1}\;\hrulefill}                  &         &         &        & \A             \\
\cline{1-7}
\cline{11-12}
    \end{array}        
    \,\right).
\end{equation}
\end{dfn}

In the above inequality system, the framed blocks denote matrices of
appropriate sizes that contain the depicted block repeated row-wise
$\frac{\gon}{2}$ times. In particular,  
\[
    \begin{array}{|cc|}
        \hline & 1 \\ & \\ \hline
    \end{array} 
    :=
    \left(\begin{array}{cc}
         0 & 1 \\
         0 & 0 \\
        \multicolumn{2}{c}{\vdots}\\
         0 & 1 \\
         0 & 0 \\
     \end{array}\right) \in \R^{\gon\times 2} 
     \ \ \text{ and } \ \ 
    \begin{array}{|cccc|}
        \hline 
        & & & 1 \\ 
\multicolumn{4}{|c|}{\hrulefill\;g_1\;\hrulefill}\\
    \hline\end{array} 
    :=
    \left(\begin{array}{cccc}
      0 & \cdots & 0 & 1 \\ 
\multicolumn{4}{c}{\hrulefill\;g_1\;\hrulefill}\\
\multicolumn{4}{c}{\vdots}\\
      0 & \cdots & 0 & 1 \\ 
\multicolumn{4}{c}{\hrulefill\;g_1\;\hrulefill}\\
     \end{array}\right) \in \R^{\gon\times (2r-d)}.
\]

\begin{proof}[Proof of Theorem \ref{thm:dpp}]
    Let $P = \dpp{2r}{\ovG}$ be the deformed product of $\gon$-gons according
    to the DPP template (\ref{eqn:dpp})  which is determined by a Gale
    transform $G = \left( \begin{smallmatrix} I_{d-2} \\ \ovG
    \end{smallmatrix} \right)$ of a neighborly $(d-2)$-polytope $Q$ with
    $2r-1$ ordered vertices in general position.  Equipped with a suitable
    right hand side, the polytope $P$ is an \emph{iterated rank 2 deformed
    product} of polygons and thus combinatorially equivalent to the $r$-fold
    product of an $\gon$-gon.

    Now for an arbitrary vertex $v$ of $P$, the matrix $\ovA_{\eq v}$ is of
    the following form
    \begin{equation}\label{eqn:G_vert}
        \ovA_{\eq v} = \left(
        \begin{array}{ccccccc}
  a_{i_1}  & a^\prime_{i_1}  &          &                 &         &                &                        \\ \hline
     1     &     a_{i_2}     &          &                 &         &                &                        \\
           &             1   & a_{i_3}  & a^\prime_{i_3}  &         &                &                        \\
           &                 &    1     &     a_{i_4}     &         &                &                        \\
           &                 &          &        \ddots   & \ddots  &                &                        \\
           &                 &          &                 &    1    & a_{i_{2r-d-1}} & a^\prime_{i_{2r-d-1}}  \\
           &                 &          &                 &         &       1        & a_{i_{2r-d}}           \\
           &                 &          &                 &         &                &        1               \\
            \multicolumn{7}{c}{\;\hrulefill\;g_1\;\hrulefill\;}                                               \\
            \multicolumn{7}{c}{\vdots}                                                                        \\
            \multicolumn{7}{c}{\;\hrulefill\;g_{d-1}\;\hrulefill\;}                                           \\
        \end{array} \right) \in \R^{2r \times (2r-d)}.
    \end{equation}
    The equality set of a vertex $v$ is formed by two cyclicly adjacent facets
    from each polygon in the product. This means, in particular, that from
    each polygon there is an even and an odd facet present in $\eq v$. Every
    such pair is of the form 
    \[
        \left(
        \begin{array}{cc}
            a_{i_\ell} & a_{i_\ell}^\prime \\
            1      & a_{i_{\ell+1}} \\
        \end{array}
        \right).
    \]
    The absolute values of the diagonal entries are bounded by $\eps$, while 
    $|a^\prime_{i_{\ell+1}}| < \eps^2$.
    
    Thus, provided that $\eps$ is sufficiently small, the rows of $\ovA_{\eq
    v}$ below the horizontal bar in (\ref{eqn:G_vert}) constitute a Gale
    transform of a polytope combinatorially equivalent to $Q$.
\end{proof}

In analogy to the cubical case, we write $\dpp{2r}{Q}$ for the deformed
product of $m$-gons with respect to the polytope $Q$ with ordered vertices.

\begin{dfn}
    The proof of Theorem \ref{thm:dpp} yields a family of \emph{projected
    products of polygons} (PDPPs) as the image $\pdpp{2r,d}{Q} :=
    \pi_d(\dpp{2r}{Q})$.
\end{dfn}

En route to a facial description of $\pdpp{2r,d}{Q}$, let us pause to
introduce a convenient notation for handling products of even polygons
combinatorially that bears certain similarities with that of $2r$-cubes, i.e.\
products of quadrilaterals. For the even polygons above, we label the edge
with outer normal $a_i$ by $(i,*)$ if $i$ is even and by $(*,i)$ otherwise:
\[
    \psfrag{P}[cb]{\scriptsize{$(*,*)$}}
    \psfrag{e1}[rb]{\scriptsize{$(0,*)$}}
    \psfrag{e2}[rb]{\scriptsize{$(*,1)$}}
    \psfrag{e3}[lb]{\scriptsize{$(2,*)$}}
    \psfrag{e4}[lb]{\scriptsize{$(*,3)$}}
    \psfrag{e5}[lt]{\scriptsize{$(4,*)$}}
    \psfrag{e6}[rt]{\scriptsize{$(*,5)$}}
    \psfrag{v12}[cr]{\scriptsize{$(0,1)$}}
    \psfrag{v23}[bc]{\scriptsize{$(2,1)$}}
    \psfrag{v34}[cl]{\scriptsize{$(2,3)$}}
    \psfrag{v45}[cl]{\scriptsize{$(4,3)$}}
    \psfrag{v56}[tc]{\scriptsize{$(4,5)$}}
    \psfrag{v16}[cr]{\scriptsize{$(0,5)$}}
    \includegraphics[height=4cm]{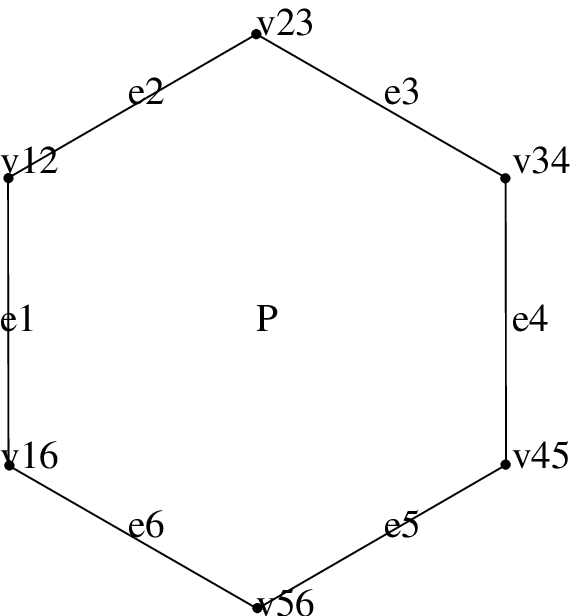}
\]
Every vertex is incident to an even edge $(2i,*)$ and an odd edge $(*,2i\pm1)$
and is labeled by $(2i,2i\pm1)$.  Finally, the polygon itself gets the label
$(*,*)$ as the intersection of no edges.

Summing up, the nonempty faces of an even $\gon$-gon are given by
\[
\begin{aligned}
    \ngonComb \ \ =\ \  & \{ (2i,*) : 0 \le i < \tfrac{\gon}{2} \} & 
                                    \text{(even edges)}\\
               \cup\ \ & \{ (*,2i+1) : 0 \leq i < \tfrac{\gon}{2}\} & 
                                    \text{ (odd edges)}\\
               \cup\ \ & \{ (2i,2i\pm1) : 0 \leq i < \tfrac{\gon}{2} \} & 
                                    \text{ (vertices)}\\
               \cup\ \ & \{ (*,*) \}  & \text{ (polygon)}\\
\end{aligned}
\]
with inclusion given by the order relation induced by $i \prec *$ for $i \in
\{0,\dots,\gon-1\}$.

Admittedly, this is neither the most natural nor the most efficient way to
encode a polygon combinatorially. However, the following remarks make up for
this unusual description. Similar to the description of $2r$-cubes, the
dimension of a face $(\a_0,\a_1) \in \ngonComb$ is the number of $*$-entries.
This carries over to products of $\gon$-gons, i.e.\ there is an
order-preserving bijection between the nonempty faces of an $r$-fold product
of $\gon$-gons and the $r$-fold direct product $(\ngonComb)^r$ with rank
function $\dim \a = \#\{ i : \a_i = * \}$ for $\a \in (\ngonComb)^r$.  Notably
most of the results (and proofs) from Section \ref{sec:ncp} carry over to this
setting, with only minor modifications. 

The key to obtaining a combinatorial description of $\pdpp{2r,d}{Q}$ is that
for a vertex $v$ of $\dpp{2r}{Q}$ the matrix (\ref{eqn:G_vert}) again encodes
a lexicographic triangulation of $Q$. In order to reduce this to the case of
neighborly cubical polytopes, after a suitable change of basis, the matrix
$\ovA_{\eq v}$ is of the form
\begin{equation}\label{eqn:ovA_form}
    \newcommand\g{\tilde{g}}
    \newcommand{\at}{\tilde{a}}
    \newcommand\ssa[1]{\scriptstyle \a_{#1}}
    \newcommand\svdots{\scriptstyle \vdots}
    \left(
    \begin{array}{ccccccc}
a_{i_1}  &                 &          &                 &         &
&                        \\ 
   1     &   \at_{i_2}     &          &                 &         &                &                        \\
         &             1   & a_{i_3}  &                 &         &                &                        \\
         &                 &    1     &   \at_{i_4}     &         &                &                        \\
         &                 &          &        \ddots   & \ddots  &                &                        \\
         &                 &          &                 &    1    & a_{i_{2r-d-1}} &                        \\
         &                 &          &                 &         &       1        & \at_{i_{2r-d}}         \\
         &                 &          &                 &         &                &        1               \\
        \multicolumn{7}{c}{\;\hrulefill\;\g_1\;\hrulefill\;}                                                \\
        \multicolumn{7}{c}{\vdots}                                                                          \\
        \multicolumn{7}{c}{\;\hrulefill\;\g_{d-1}\;\hrulefill\;}                                            \\
    \end{array} \right)
    \begin{array}{l}
        \ssa{1} \\
        \ssa{2} \\
        \ssa{3} \\
        \ssa{4} \\
        \svdots \\
        \ssa{2r-d-1} \\
        \ssa{2r-d} \\
        \ssa{2r-d+1} \\
        \ssa{2r-d+2} \\
        \svdots \\
        \ssa{2r} \\
    \end{array}
\end{equation}
The entries above the diagonal of ones remain to be of order $\eps$.  To
determine the signs of the entries, which will determine the lexicographic
triangulation, let us investigate the \emph{local} change of the matrix under
the change of basis.

In the above combinatorial model for even $\gon$-gons, the vertex $v$ is
identified with a vector $\a =  (\a_1,\a_2;\a_3, \dots;\a_{2r-1},\a_{2r}) \in
(\ngonComb)^r$, which corresponds to $\eq v$ as indicated. The following
table, which is easily established given the coordinates of the normals,
summarizes the possible sign patterns in terms of $\a$.
\[
    \begin{array}{c||c|c|c|c}
        (\a_i,\a_{i+1}) & (0,1) & (0,\gon-1) & (2k,2k-1) & (2k,2k+1) \\
        \hline
        \left( \begin{array}{cc}
           a_i   &        \\
           1 & \tilde{a}_{i+1} \\
        \end{array} \right)
        &
        \left( \begin{array}{cc}
           -\eps &        \\
               1 & +\eps   \\
        \end{array} \right)
        &
        \left( \begin{array}{cc}
           -\eps &        \\
               1 & -\eps  \\
        \end{array} \right)
        &
        \left( \begin{array}{cc}
            +\eps &      \\
               1 & +\frac{2\eps}{\gon-2}   \\
        \end{array} \right)
        &
        \left( \begin{array}{cc}
            +\eps &      \\
               1 & -\frac{2\eps}{\gon-2}   \\
        \end{array} \right) \\
        \hline
            (\s_i,\s_{i+1}) & (-,+) & (-,-) & (+,+) & (+,-) \\
    \end{array}
\]

We use the last row, which gathers sign patterns from the diagonal, to define
the map 
\[
    \Phi: \{ (\a_1,\a_2) \in \ngonComb : \a \text{ vertex} \} \rightarrow \PMZ^2
\]
with $\Phi(\a_1,\a_2) := (\s_1,\s_2)$ according to the table. Since the face
lattice of a convex polytope is atomic, it is easy to see from the definition
that $\Phi : \ngonComb \rightarrow \PMZ^2$ extends to an order- and
rank-preserving map from the face poset of an even $\gon$-gon to that of an
$2$-cube. The map can be thought of as a \emph{folding} map: 
\[
    \psfrag{Phi}[lc]{$\Phi$}
    \psfrag{++}[lc]{$\scriptstyle ++$}
    \psfrag{+-}[lc]{$\scriptstyle +-$}
    \psfrag{--}[rc]{$\scriptstyle --$}
    \psfrag{-+}[rc]{$\scriptstyle -+$}
    \psfrag{05}[rb]{$\scriptstyle (0,5)$}
    \psfrag{01}[rc]{$\scriptstyle (0,1)$}
    \psfrag{21}[cb]{$\scriptstyle (2,1)$}
    \psfrag{23}[lb]{$\scriptstyle (2,3)$}
    \psfrag{43}[lb]{$\scriptstyle (4,3)$}
    \psfrag{45}[cc]{$\scriptstyle (4,5)$}
    \includegraphics[width=12cm]{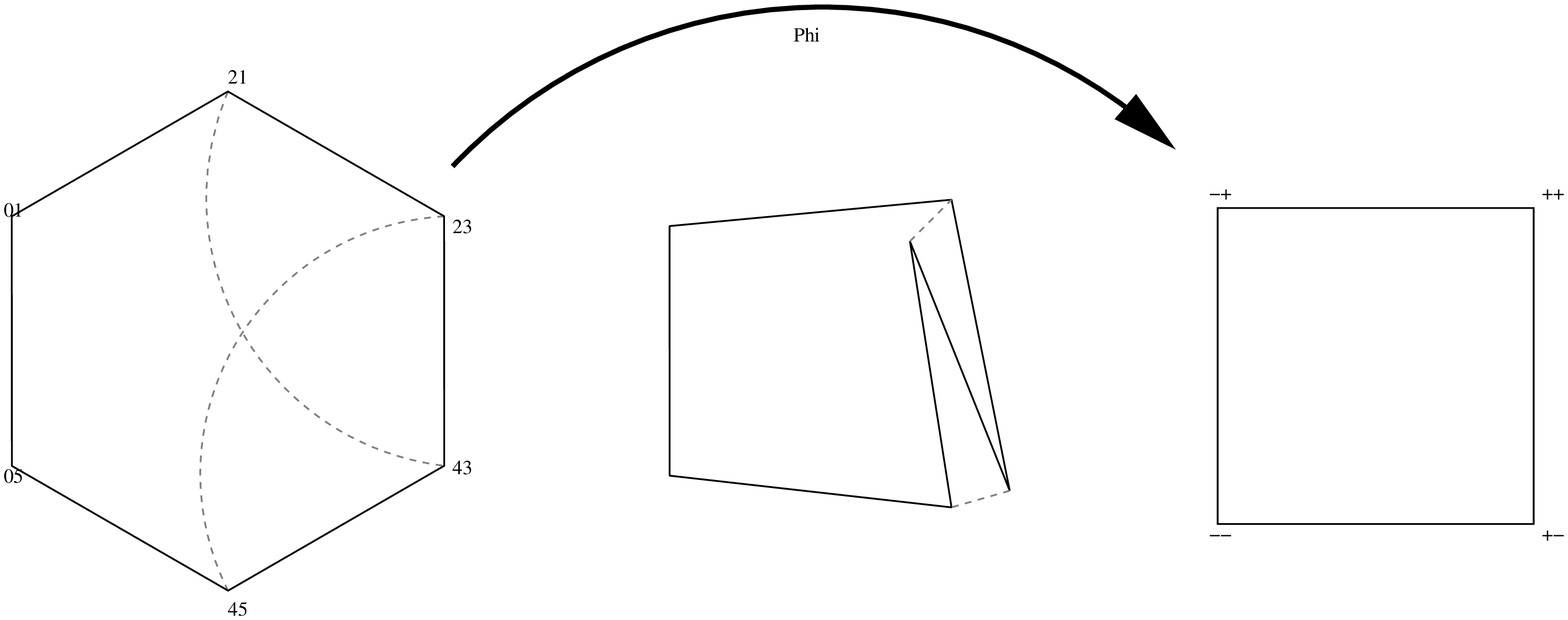}
\]

The induced map $\Phi : (\ngonComb)^r \rightarrow \PMZ^{2r}$ maps faces of
$\dpp{2r}{k}$ that are strictly preserved under $\pi_d$ to surviving faces of
$C_{2r}(Q)$. Phrased differently the following diagram commutes on the level
of faces:
\[
\begin{CD}
    P_{n,r}(Q) @>{\hspace{1cm}\Phi\hspace{1cm}}>> C_{2r}(Q) \\
      @VV\pi_dV                                                 @VV\pi_dV   \\
    \pdpp{2r,d}{Q} @>{\hspace{1cm}\Phi\hspace{1cm}}>> \ncp_d(Q).
\end{CD}
\]
\begin{prop}
    Let $n = 2r$ and let $P = \dpp{n}{Q}$ and $C = C_{n}(Q)$ be the deformed
    cube and the product of $\gon$-gons of dimension $n = 2r$ with respect to
    a neighborly $(d-2)$-polytope $Q$ on $n-1$ ordered vertices.  Let $v \in
    P$ be a vertex with $\eq v$ represented by $\a \in (\ngonComb)^r$ and let
    $u \in C$ be the vertex corresponding to $\Phi(\a) \in \PM^{n}$. Then
    $\Phi$ induces an isomorphism of the vertex figures $\pdpp{n,d}{Q} /
    \pi_d(v)$ and $\ncp_n(Q) / \pi_d(u)$.
\end{prop}
\begin{proof}
    As consistent with the main theme in this article, consider the first $n-d
    = 2r-d$ coordinates of the inequalities from both $P$ and $C$ that are
    tight at $v$ and $u$, respectively.
\[
    \newcommand\se[1]{\s_{i_{#1}}\eps}
    \newcommand\g{\tilde{g}}
    \newcommand{\at}{\tilde{a}}
    \newcommand\ssa[1]{\scriptstyle \a_{#1}}
    \newcommand\svdots{\scriptstyle \vdots}
    \begin{array}{c@{\hspace{1cm}}c}
    \ovA_v(P) &
    \ovA_u(C) \\
    \left(
    \begin{array}{ccccc}
a_{i_1}  &                  &                 &         &              \\ \hline
   1     &   \at_{i_2}      &                 &         &              \\
         &         \ddots   & \ddots  &                &               \\
         &                  &    1    & a_{i_{n-d-1}} &               \\
         &                  &         &       1        & \at_{i_{n-d}}\\
         &                  &         &                &        1      \\
        \multicolumn{5}{c}{\;\hrulefill\;\g_1\;\hrulefill\;}           \\
        \multicolumn{5}{c}{\vdots}                                     \\
        \multicolumn{5}{c}{\;\hrulefill\;\g_{d-1}\;\hrulefill\;}       \\
    \end{array} \right)
    &
    \left(
    \begin{array}{ccccc}
\se{1}   &             &         &             &           \\ \hline
   1     &   \se{2}    &         &             &           \\
         &    \ddots   & \ddots  &             &           \\
         &             &    1    & \se{n-d-1} &           \\
         &             &         &       1     & \se{n-d} \\
         &             &         &             &        1  \\
        \multicolumn{5}{c}{\;\hrulefill\;g_1\;\hrulefill\;}            \\
        \multicolumn{5}{c}{\vdots}                                     \\
        \multicolumn{5}{c}{\;\hrulefill\;g_{d-1}\;\hrulefill\;}        \\
    \end{array} \right)
    \end{array} 
\]
In both matrices, the entries on the secondary diagonal are arbitrary small
and the map $\Phi$ assures that corresponding entries have equal sign. By
Lemma \ref{lem:lex_gale}, both $\ovA_v(P)$ and $\ovA_u(C)$ are Gale transforms
that encode the same lexicographic pyramid $\L_k(Q)$. The result now follows
by observing that a face $\b \succeq \a$ of $P$ is strictly preserved if and
only if $|\b|$ is a coface of $\L_k(Q)$ and $|\Phi(\b)| = |\b|$.
\end{proof}

This proposition makes way for the combinatorics of the projected deformed
products associated with arbitrary simplicial neighborly polytopes.

\begin{thm}[Combinatorial Description of the PDPPs]\label{thm:pdpp_comb}
    Let $P = \dpp{2r}{Q}$ be a deformed product of $\gon$-gons with respect to
    $Q$ and let $v \in P$ be an arbitrary vertex with $\eq v = \a \in
    (\ngonComb)^r$. Then the vertex figure of $\pi_d(v)$ in $\pdpp{2r,d}{Q}$
    is isomorphic to $\L_p(Q)$ for
    \[
        p = \min\{ i \in [2r] : \Phi(\a)_i = - \} \cup \{ 2r-d-1 \}.
    \]
    In particular, the $(d-1)$-faces of $P$ containing $v$ that are preserved
    by projection are in one-to-one correspondence to the facets of $\L_p(Q)$.
\end{thm}

As for the neighborly cubical polytopes, via Shemer's work \cite{shemer} this
result implies a great richness of combinatorial types for the projected
products of polygons.  In the special case when $Q$ is a cyclic polytope with
vertices in standard order, we get a very explicit \emph{Gale's evenness}-type
criterion for the projected products of polygons.

\begin{cor}[Combinatorial Description of the standard PDPPs]
    Let $F \subset P = \dpp{2r}{Q}$ be a $(d-1)$-face with $Q =
    \cyc_{d-2}(2r-1)$ and let $\b \in (\ngonComb)^r$ correspond to $\eq F$.
    Then $F$ projects to a facet of $\pdpp{2r,d}{Q}$ if and only if $\Phi(\b)$
    satisfies the cubical Gale's evenness criterion.
\end{cor}

\begin{small}

\end{small}

\end{document}